\documentclass[graybox]{svmult}
\usepackage{newtxtext,newtxmath,amssymb}
\usepackage{siunitx}
\usepackage{helvet}         
\usepackage{courier}        
\usepackage{makeidx}         
\usepackage{graphicx}        
\usepackage[bottom]{footmisc}
\usepackage{url}
\makeindex

\usepackage{amsfonts}%
\usepackage{float}
\usepackage{bbold}
\usepackage{color}
\usepackage{caption}
\usepackage{url}

\theoremstyle{plain}

\begin{document}

\title*{\textbf{ A survey on knotoids, braidoids and their applications}}
 
\author{Neslihan G\"ug\"umc\"u, Louis H.Kauffman and Sofia Lambropoulou}
\institute{Neslihan G\"ug\"umc\"u \at {School of Applied Mathematical and Physical Sciences \\ National Technical University of Athens,
Zografou Campus, GR-157 80 Athens, Greece. Greece} \\ \email{nesli@central.ntua.gr} \quad url: \url{www.nesligugumcu.weebly.com}
\and Louis H.Kauffman  \at {Department of Mathematics, Statistics, and Computer Science, University of Illinois at Chicago, 851 South Morgan St., Chicago IL 60607-7045,  U.S.A. } \\ Department of Mechanics and Mathematics,  Novosibirsk State University, Novosibirsk, Russia. \\ \email{kauffman@uic.edu} \quad url:\url{http://homepages.math.uic.edu/~kauffman/}
\and Sofia Lambropoulou  \at School of Applied Mathematical and Physical Sciences \\ National Technical University of Athens, Zografou Campus, GR-157 80 Athens, Greece. \\ \email{sofia@math.ntua.gr} \quad url: \url{http://www.math.ntua.gr/~sofia}}
\maketitle

 





\let\thefootnote\relax\footnotetext{

{\noindent}\textit{2010 Mathematics Subject Classification}: 57M27, 57M25 \\
{\noindent}\textit{Keywords}: knotoids, multi-knotoids, spherical knotoids, planar knotoids, loop bracket polynomial, arrow loop polynomial,  $\Theta$-graphs, line isotopy, braidoiding algorithm, $L$-braidoiding move, Alexander theorem, $L$-moves, Markov theorem, proteins}

\abstract{
This paper is a survey on the theory of knotoids and braidoids. Knotoids are open ended knot diagrams in surfaces and braidoids are geometric objects analogous to classical braids, forming a counterpart theory to the theory of knotoids in the plane. We survey through the fundamental notions and existing works on these objects  as well as their applications in the study of proteins.
}
\setcounter{section}{-1} 

\section{Introduction}\label{intro}

The theory of knotoids was introduced by Turaev \cite{Tu} in 2010. A  knotoid diagram is an oriented curve with two endpoints, immersed in an oriented surface, and having finitely many self-intersections that are endowed with under/over data. The two endpoints may lie in different regions of the diagram. They may move within their regions by planar isotopy, but they are not allowed to cross over or under any arcs of the diagram. These are the `forbidden moves' of the theory. The  definition of a  knotoid diagram  extends the notion of $(1, 1)$-tangle whose endpoints can be assumed to be fixed at the boundary of the disk where the tangle lies. The theory of  knotoids in the 2-sphere extends the theory of classical knots \cite{Tu} and also proposes a new diagrammatic approach to  classical knot theory and  to  classical knot invariants, with reduced computational complexity, at least for invariants computed on the number of crossings of a given diagram, see  \cite{Tu,DST}. In \cite{Tu} basic properties of knotoids were studied, including the introduction of several invariants of spherical knotoids,  such as the complexity (or height, as renamed in \cite{GK1}) and the Jones/bracket polynomial. Knotoids in $S^2$ up to 5 crossings  were classified by Bartholomew in 2010, using the generalization of the bracket polynomial for knotoids that was defined by Turaev  \cite{Ba}.  In 2013 the subject of knotoids and the introduction of the theory of braidoids was  proposed to the first listed author by the third one for her PhD study.  Then, the first and the second listed authors introduced several new invariants  in analogy with invariants from virtual knot theory \cite{Gthesis,GK1}, see also remark in \cite{Tu}, while in \cite{Gthesis,GL1,GL2}  the theory of braidoids was initiated by the first and the last authors in relation to the theory of planar knotoids. 

 A braidoid diagram extends the notion of classical braid diagram \cite{Ar1,Ar2} with extra `free' strands that initiate/terminate at two endpoints located anywhere in the plane of the diagram. The  closure operation for braidoids requires special attention due to the presence of the endpoints and their forbidden moves, while a `braidoiding' algorithm,  based on earlier work of the third listed author \cite{Lathesis,LR1}, that turns any planar knotoid diagram into a braidoid diagram is the proof of an analogue of the classical Alexander theorem~\cite{Al} for knotoids. Further, with the introduction of $L$-moves on braidoids, which were originally defined for classical braids by the last author \cite{Lathesis,LR1}, a geometric analogue of the classical Markov theorem~\cite{Ma} for braidoids is enabled. In \cite{Gthesis,GL1,Gu1} a set of combinatorial elementary blocks for braidoids is also introduced, which in \cite{GL1}  are proposed for an algebraic encoding of the entanglement of open protein chains in $3$-dimensional space.  

 Recently, knotoids have been used in the field of biochemistry as they suggest new topological models for open linear molecules. More precisely, the invariants for spherical and planar knotoids introduced in \cite{Tu, Gthesis, GK1} have been used for determining the topological entanglement of open protein chains in \cite{GDBS, GGLDSK,DRGDSMRSS}.   In particular, in \cite{GGLDSK} the concept of `bonded knotoid' has been introduced for this purpose. Some other recent works on knotoids include a study of biquandle coloring invariants, by the first listed author and Nelson \cite{GN1}; the study of knots that are knotoid equivalent, by Adams, Henrich, Kearney and Scoville \cite{AdHen};   a recent classification table for prime knotoids of positive complexity with up to 5 crossings \cite{KMT} given by Korablev, May and Tarkaev, obtained by using the correspondence between the knotoids in $S^2$ and the knots in thickened torus; a classification of all planar knotoids with up to 5 crossings and extension of spherical-knotoids to 6 crossings; the introduction of the equivalent theory of `rail knotoids' and their connection to the knot theory of the handlebody of genus~2 by Kodokostas and the third listed author \cite{KoLa1,KoLa2}; and the construction of double branched covers of knotoids \cite{BBHL}. 


The outline of the paper is as follows. In Section~\ref{sec:knotoids} we review the basic notions of the theory of knotoids. In Section~\ref{sec:closures} we present closure types for knotoids for obtaining knots. In Section~\ref{sec:sph} we focus on the spherical knotoids and how they extend the classical knot theory. In Section~\ref{sec:geo} we present geometric interpretations for  spherical and planar knotoids. In Section~\ref{sec:inv} we survey through the existing works and results on the invariants of knotoids. In Section~\ref{sec:braidoids} we review the  notion of braidoid. In Sections~\ref{sec:closure} and~\ref{sec:braiding}, we present the key elements for proving an analogue of the Alexander theorem for knotoids/multi-knotoids. In Section~\ref{sec:l} we present the definition of the $L$-moves on braidoid diagrams that give rise to an analogue of the Markov theorem for braidoids. Finally, in Section~\ref{sec:app} we present the applications of knotoids to the study of proteins; we also review the building blocks for braidoid diagrams,  proposed to be used in encoding open protein chains.

\section{Knotoids and knotoid isotopy} \label{sec:knotoids}

\subsection{Knotoid diagrams}

\begin{figure}[H]
\centering  \scalebox{0.8}{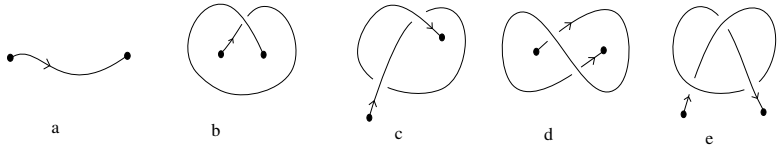}
\caption{Examples of knotoid diagrams}
\label{fig:knotoid}
\end{figure}

Let $\Sigma$ be an oriented surface. A \textit{knotoid diagram} $K$ in $\Sigma$ \cite{Tu} is an immersion of the unit interval $[0,1]$ in $\Sigma$ with a finite number of double points each of which is a transversal self-intersection endowed with over/under data. Such self-intersections of $K$ are called {\it crossings} of $K$. The images of $0$ and $1$ are two distinct points called the \textit{endpoints} of $K$ and are specifically called the \textit{leg} and the \textit{head}, respectively. A knotoid diagram is naturally oriented from its leg to its head. The trivial knotoid diagram is assumed to be an immersed curve without any self-intersections as depicted 
in Figure~\ref{fig:knotoid}a. 
 
The notion of knotoid can be extended to include more components. A \textit{multi-knotoid diagram} in $\Sigma$ is a union of a knotoid diagram and a finite number of knot diagrams \cite{Tu}. 


\subsection{Moves on knotoid diagrams} 
\begin{figure}[H]
 \centering  \scalebox{0.45}{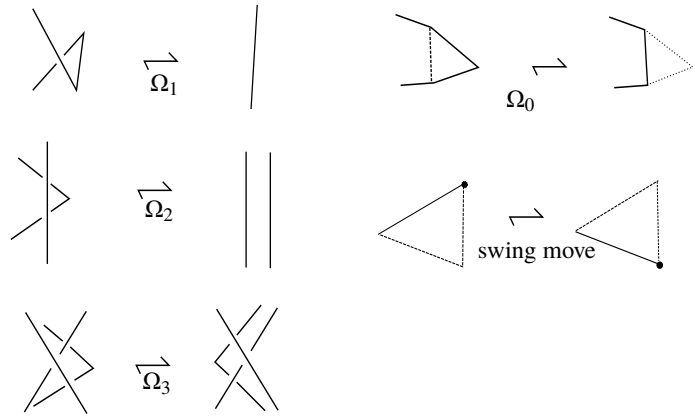}
\caption{The moves generating isotopy on knotoid diagrams}
       \label{fig:om}
	\end{figure}

Planar isotopy moves generated by the $\Omega_0$-move and the Reidemeister moves $\Omega_1$, $\Omega_2$, $\Omega_3$ (see Figure~\ref{fig:om}) that take place in a local disk free of any endpoints are allowed on knotoid diagrams. A special case of planar isotopy moves that involves an endpoint is a \textit{swing move}, whereby an endpoint can be pulled within its region, without crossing any other arc of the diagram, as illustrated in Figure~\ref{fig:om}. 
We refer to all these moves as $\Omega$-moves of knotoids.

\begin{figure}[H]
       \centering  \scalebox{0.6}{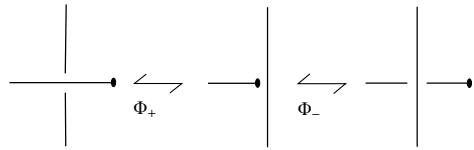}
				\caption{  Forbidden knotoid moves}
        \label{subfig:for}
		\end{figure}

The moves consisting of pulling the arc adjacent to an endpoint over or under a transversal arc, as shown in Figure~\ref{subfig:for}, are the \textit{forbidden knotoid moves} and are denoted by $\Phi_+$ and $\Phi_-$, respectively. Notice that, if both $\Phi_+$ and $\Phi_-$-moves were allowed, any knotoid diagram in any surface could be clearly turned into the trivial knotoid diagram. 

\subsection{Knotoids} \label{sec:knotoid}
  The  $\Omega$-moves generate an equivalence relation on knotoid diagrams in $\Sigma$, called {\it knotoid isotopy}. Two knotoid diagrams are {\it isotopic} to each other if there is a finite sequence of $\Omega$-moves that takes one to the other. The isotopy classes of knotoid diagrams in $\Sigma$ are called \textit{knotoids}. The equivalence relation defined on knotoid diagrams applies also to multi-knotoid diagrams, and the corresponding equivalence classes are called \textit{multi-knotoids}. 

Let $\mathcal{K}(\mathbb{R}^2)$ and $\mathcal{K}(S^2)$ denote the set of all knotoids in $\mathbb{R}^2$ and $S^2$, respectively. We shall call knotoids in $\mathcal{K}(\mathbb{R}^2)$ {\it planar} and knotoids in $\mathcal{K}(S^2)$ {\it spherical}.  

 There is a surjective map  \cite{Tu}
$$
\iota:\mathcal{K}(\mathbb{R}^2)\hookrightarrow \mathcal{K}(S^2),
$$ 
induced by the inclusion $\mathbb{R}^2\hookrightarrow S^2=\mathbb{R}^2\cup\{\infty\}$. However, the map $\iota$ is not injective.  Indeed, there are knotoid diagrams in the plane representing a nontrivial knotoid while they represent the trivial knotoid in $S^2$. For an example, see Figure~1(b). It is well-known that the knot theory of the plane coincides with the knot theory of the $2$-sphere, while the non-injectivity of the map $\iota$ implies that the theory of knotoids in $\mathbb{R}^2$ differs from the theory of knotoids in $S^2$, also yields a more refined theory than the theory of spherical knotoids.

\section{Knotoids, classical knots and virtual knots} \label{sec:closures}

\subsection{Classical knots via knotoids}
\begin{figure}[H]
\centering  \scalebox{0.27}{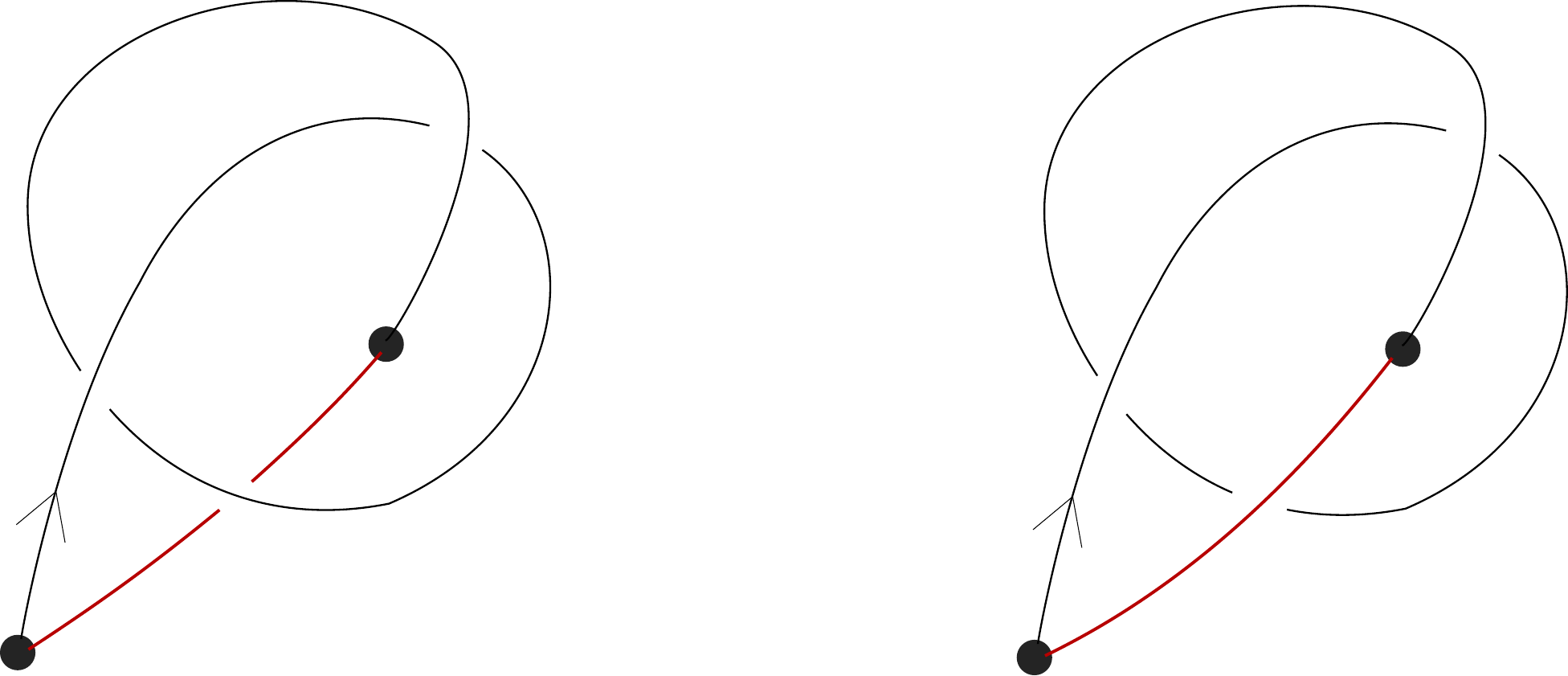}
\caption{ The overpass and the underpass closures of a knotoid diagram resulting in different knots}
\label{fig:cl}
\end{figure}

In \cite{Tu} the study of knotoid diagrams is suggested as a new diagrammatic approach to the study of knots in 3-dimensional space, as any classical knot can be represented by a knotoid diagram in $\mathbb{R}^2$ or in $S^2$. More precisely, the endpoints of a knotoid diagram can be connected with an arc in $S^2$ that goes either under each arc it meets or over each arc it meets, as illustrated in  Figure~\ref{fig:cl}. This way we obtain an oriented classical knot diagram in $S^2$ representing a knot in $\mathbb{R}^3$. The connection types are called the \textit{underpass} \textit{closure} and the\textit{ overpass closure}, respectively. The knot that is represented by a knotoid diagram may differ depending on the type of the closure. For example, the knotoid in Figure~\ref{fig:cl} represents a trefoil via the underpass closure and  the trivial knot via the overpass closure. 

  In order to have a well-defined representation of knots via knotoids, we should fix the closure type for knotoid diagrams. When we choose the underpass closure as closure type, we have the following proposition. The statement of the proposition is symmetric for the overpass closure.  
	
\begin{proposition}[\cite{Tu}]
Two knotoid diagrams in $S^2$ or $\mathbb{R}^2$ represent the same classical knot if and only if they are related to each other by finitely many $\Omega$-moves, swing moves and the forbidden $\Phi_-$-moves.
\end{proposition} 

Given a knot in $S^3$. We can also consider cutting out an underpassing or an overpassing strand of an oriented diagram of the knot. The resulting diagram is clearly a knotoid diagram in the plane or the $2$-sphere. In fact we can obtain a set of knotoid diagrams by cutting out different strands and the resulting knotoid diagrams all represent the given knot via the underpass or the overpass closure. Knotoid representatives of a knot in $S^3$ clearly have less number of crossings than any of the diagrams of the knot. For this reason, use of knotoid diagrams to study knots in $S^3$ provides a considerable amount of simplification for computing knot invariants, such as the knot group \cite{Tu}. 

Another interesting question in relation with the knotoid closures has been recently worked in \cite{AdHen}. Two knots $K_1, K_2$ in $S^3$ are said to be \textit{knotoid equivalent} if there exists a knotoid $\kappa$ such that $K_1$ is the underpass closure of  $\kappa$ and $K_2$ is the overpass closure of $\kappa$. So the question is: \textit{Which pairs of knots are knotoid equivalent?}  The authors proved the following theorem.
 
\begin{theorem}[\cite{AdHen}]
Given any two knots $K_1, K_2$ in $S^3$, $K_1$ is knotoid equivalent to $K_2$.
\end{theorem}

\subsection{Virtual knots via knotoids}

\begin{figure}[H]
\centering
\includegraphics[width=.6\textwidth]{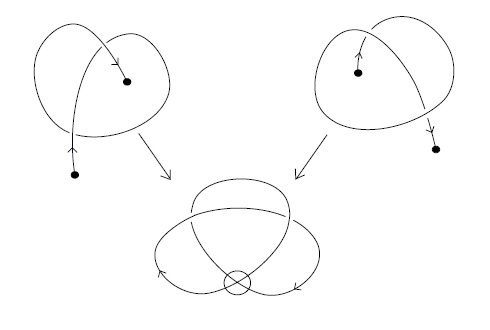}
\caption{ The virtual closure of two knotoids}
\label{fig:v}
\end{figure}

The endpoints of a knotoid diagram in $S^2$ can be tied up also in the virtual way (observation first made by Viro as mentioned in \cite{Tu}). Namely, the endpoints of a knotoid diagram can be connected with an arc by declaring each intersection of the arc with the diagram as a virtual crossing, as illustrated in Figure~\ref{fig:v}. This induces a non-injective (e.g. the knotoids in Figure~\ref{fig:v} are different) and non-surjective  mapping from the set of knotoids in $S^2$ to the set of virtual knots of genus $1$, called the \textit{virtual closure map}. Being a well-defined map, the virtual closure map provides a way to extract invariants for knotoids from the virtual knot invariants \cite{Tu,GK1}. In fact, most of the new invariants of knotoids constructed in \cite{GK1} and are briefly mentioned in Section~\ref{sec:sph}, are the result of using the principle that the virtual knot class of the virtual closure of a knotoid is an invariant of the knotoid.
 
\subsection{Spherical knotoids extend the classical knot theory} \label{sec:sph}

There is a well-defined injective map from the set of oriented classical knots to $\mathcal{K}(S^2)$ \cite{Tu}. This map is induced by specifying an oriented diagram for a given oriented knot in $S^3$ and cutting out an open arc from this diagram that does not contain any crossings. The resulting diagram is a knotoid diagram in $S^2$ with its endpoints lying in the same local region of $S^2$. Such a knotoid diagram is a \textit{knot-type knotoid diagram} and the isotopy class of the diagram is a {\it knot-type knotoid}. Figures~1(a),~1(b) and~1(e) are some examples of knot-type knotoid diagrams. It is clear that this map gives a one-to-one correspondence between the set of oriented classical knots and the set of knot-type knotoids in $S^2$. 
There are also knotoids that do not lie in the image of this map. They are called \textit{proper knotoids}. The endpoints of any of the representative diagram of a proper knotoid can lie in any but different local regions of the diagram. Figures~1(c) and~1(d) illustrate some examples of proper knotoids.

\subsubsection{The monoid of knotoids}

As Turaev explains in \cite{Tu}, two knotoids $K_1$, $K_2$ in $S^2$ can be concatenated end-to-end in the following way. One can cut out regular disk-neighborhoods of the head of $K_1$ and the leg of $K_2$ and identify the remaining surfaces with boundary along their boundaries with an orientation-reversing homeomorphism  carrying the unique intersection point of $K_1$ with the regular neighborhood of the head of $K_1$ to the unique intersection point of $K_2$ with the regular neighborhood of the leg of $K_2$. The resulting diagram is a composite knotoid diagram in $S^2$, denoted by $K_1 \# K_2$. Equipped with the binary operation $\#$, the set of spherical knotoids, $K(S^2)$, carries a monoid structure \cite{Tu}.      

\section{Geometric interpretations of knotoids}\label{sec:geo}

\subsection{A geometric interpretation of spherical knotoids}

\begin{figure}[H]
\centering
\includegraphics[width=.75\textwidth]{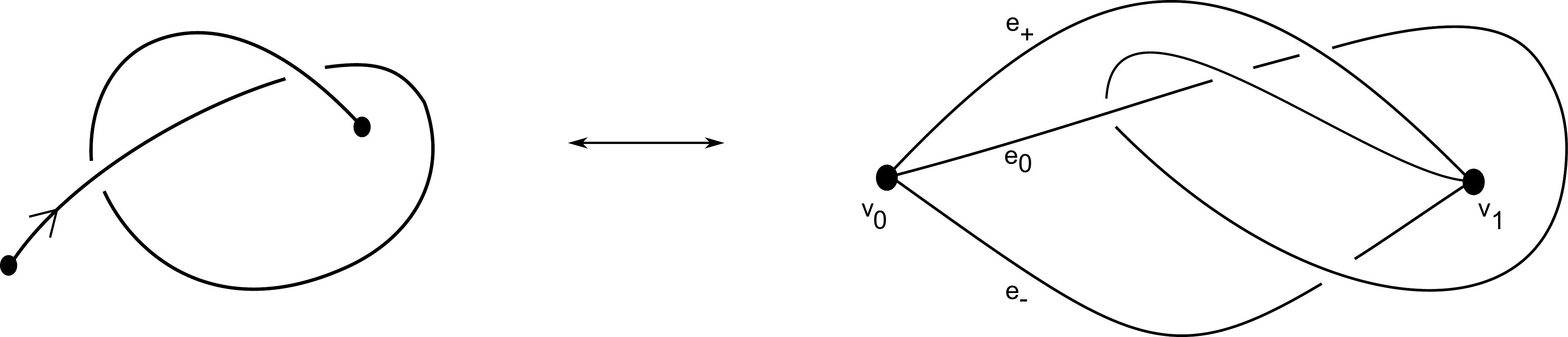}
\caption{A knotoid and the corresponding simple $\Theta$-graph}
\label{fig:theta}
\end{figure}

A {\it $\Theta$-graph} is a spatial graph with two vertices $v_0, v_1$, called the \textit{leg} and the \textit{head} respectively, and three edges, $e_+, e_0, e_-$ connecting $v_0$ to $v_1$, as exemplified in Figure~\ref{fig:theta}. The isotopy on $\Theta$-graphs is defined to be the ambient isotopy of the 3-dimensional space preserving the labeling of vertices and the edges, and the set of $\Theta$-graphs consists of the isotopy classes of $\Theta$-graphs.

 There is a binary operation on the set of $\Theta$-graphs, called the \textit{vertex multiplication} given as follows  \cite{Wo,Tu}. Let $\Theta_1$ and $\Theta_2$ be two $\Theta$-graphs, take out an open 3-disk neighborhood of the head of $\Theta_1$ and an open 3-disk neighborhood of the leg of $\Theta_2$, each intersecting with the graphs along 3-radii (simple parts from $e_0, e_+, e_-$).  Then identify the remaining manifolds with boundary along their boundaries with an orientation-reversing homeomorphism. With the vertex multiplication, the set of $\Theta$-graphs forms a monoid.
A $\Theta$-graph is called {\it simple} $\Theta$-graph if the union of its edges, $e_+$ and $e_-$ (the upper and the lower edge, respectively) is the trivial knot. Simple $\Theta$-graphs form a submonoid \cite{Tu}.

Turaev showed the following correspondence between the spherical knotoids and the simple $\Theta$-graphs that gives also rise  to a geometric interpretation of spherical knotoids via $\Theta$-curves.

\begin{theorem}[\cite{Tu}]
There is an isomorphism of monoids of spherical knotoids and of simple $\Theta$-graphs.
\end{theorem}

\subsection{A geometric interpretation of planar knotoids}

\begin{figure}[H]

\centering
\scalebox{.25}{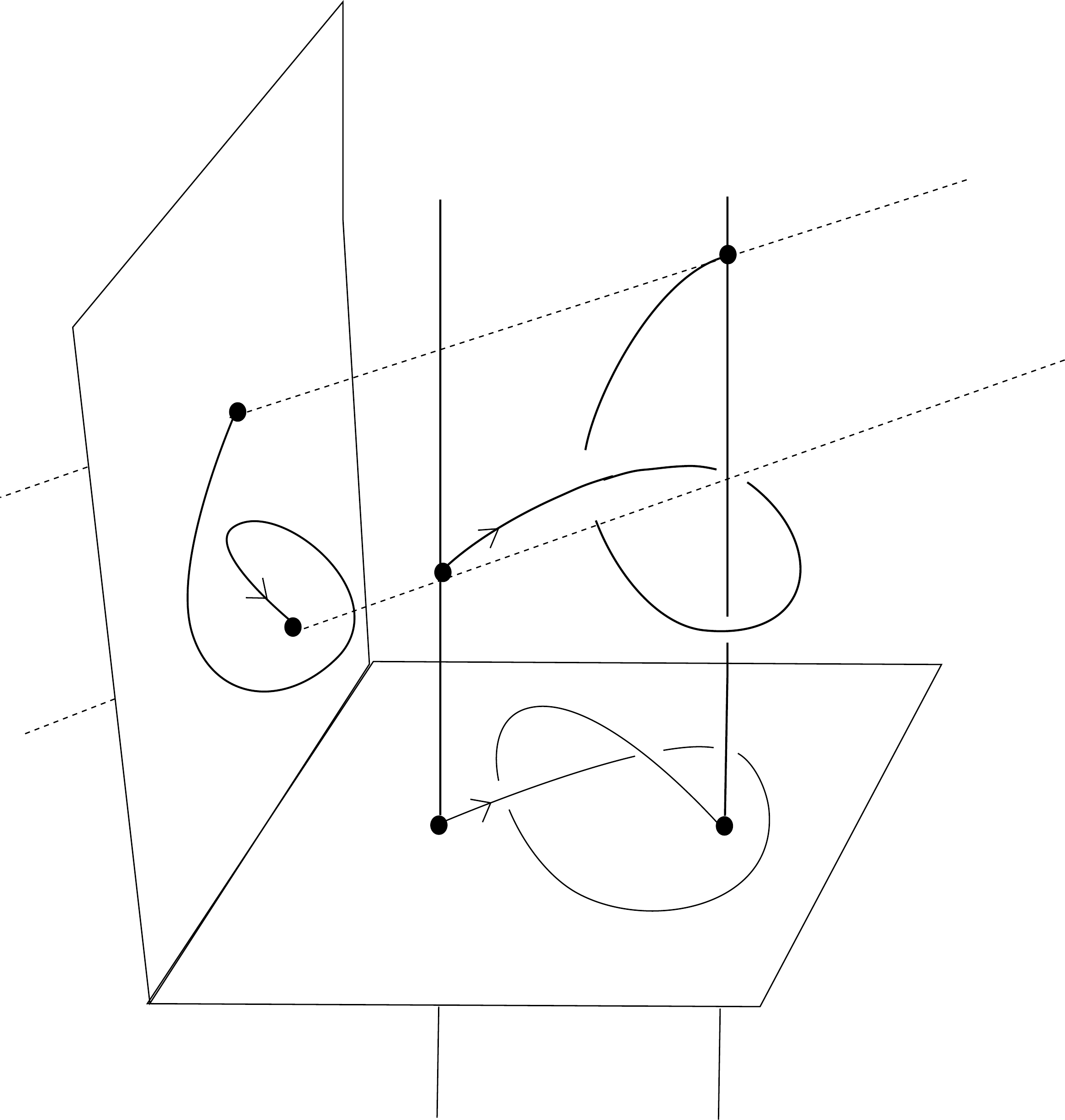}
\caption{ Projections of an open space curve as a knotoid diagram }
\label{projlines}
\end{figure}

It is explained in \cite{GK1} that the theory of planar knotoids is related naturally to open space curves on which an appropriate isotopy is defined. An open curve located in $\mathbb{R}^3$ corresponds to a planar knotoid diagram when projected regularly along the two lines passing through its endpoints and are perpendicular to a chosen projection plane. View Figure~\ref{projlines}. Conversely, an open space curve with two specified parallel lines passing from its endpoints can be viewed as a {\it lifting} of the related knotoid diagram.
 A {\it line isotopy} \cite{GK1}  between two open space curves is an ambient isotopy of $\mathbb{R}^3$ transforming one curve to the other one in the complement of the lines, keeping the endpoints on the lines and fixing the  lines. The isotopy classes of planar knotoids (considered in the chosen projection plane)  are in one-to-one correspondence with the line isotopy classes of open space curves \cite{GK1}. Furthermore, in \cite{KoLa1,KoLa2} Kodokostas and the third author make the observation that 
 this interpretation of planar knotoids as space curves is related to the knot theory of the handlebody of genus two and they propose the construction of knotoid invariants through this connection. 

\section{Invariants of knotoids}\label{sec:inv}
\begin{figure}[H]
\begin{center}
     \begin{tabular}{c}
		\Large{
     \centering  \scalebox{0.6}{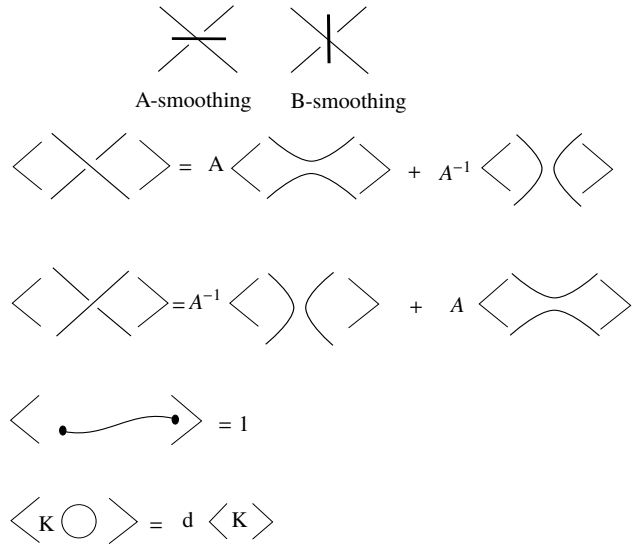}
		}
     \end{tabular}
     \caption{\bf  State expansion of the bracket polynomial}
     \label{fig:bracket}
\end{center}
\end{figure}

In \cite{Tu} and \cite{GK1} several invariants for knotoids are introduced. One of the first invariants introduced by Turaev \cite{Tu} was the the {\it bracket} and the {\it Jones polynomial} for knotoids. The bracket polynomial extends to spherical knotoids in the natural way. More precisely, the bracket expansion is directly applied to knotoid diagrams as shown in Figure~\ref{fig:bracket},  and in each state we observe a single long segment component with endpoints and a finite number of circular components. Each  circular component contributes to the polynomial by the value $-A^{-2}-A^{2}$. The initial conditions given in Figure~\ref{fig:bracket} are sufficient for the computation of the bracket polynomial of a knotoid. The closed formula for the bracket polynomial of knotoids is as follows.

\begin{definition}\normalfont
The {\it bracket polynomial} of a knotoid diagram $K$ is defined as
\begin{center}
$<K>=\sum_S A^{\sigma(S)}d^{\|S\|-1}$,
\end{center}
where the sum is taken over all states, $\sigma(S)$ is the sum of the labels of the state $S$, $\|S\|$ is the number of components of $S$, and $d=-A^2-A^{-2}$.
\end{definition}

 The bracket polynomial of knotoids in $S^2$ normalized by the writhe factor, $(-A^3)^{-\mathrm{wr}(K)}$, generalizes the Jones polynomial with the substitution $A=t^{-{1/4}}$. Note that the Jones polynomial of the trivial knotoid is trivial. Furthermore, the following conjecture \cite{GK1} extends the long-standing Jones polynomial conjecture.

\begin{conjecture}
The Jones polynomial of spherical knotoids detects the trivial knotoid.
\end{conjecture} 

Some other generalizations of the Jones polynomial are: Turaev's \textit{ 2-variable bracket polynomial} \cite{Tu} that is obtained by a use of the intersection number of the connection arc used in the underpass closure with the rest of the diagram and also with the state components, and the {\it arrow polynomial} \cite{GK1} that keeps track of the cusp-like structure (see Figure~\ref{fig:cusp}) arising in the oriented bracket expansion (see Figure~\ref{fig:arrs} for the oriented state expansion) by assigning new variables namely $\Lambda_i$'s to zig-zag components. There is a special generalization of the bracket polynomial for planar knotoids induced by distinguishing the circular state components nested around the long segment component from the circular state components that do not nest around the long segment component. Turaev defined the \textit{$3$-variable bracket polynomial} \cite{Tu} for planar knotoids based on this idea and the \textit{loop bracket polynomial}, which is a specialization of the 3-variable bracket polynomial, was utilized in \cite{GGLDSK} to classify the knotoid models of protein chains (see Section~\ref{sec:app}). Similarly the first and the second listed authors introduce the \textit{arrow loop polynomial} in \cite{GK1}.  

\begin{figure}[H]
\centering
\scalebox{.3}{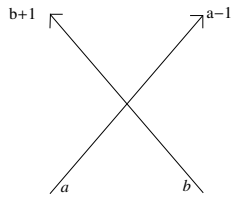}
\caption{An affine biquandle coloring on knotoid diagrams}
\label{fig:aff}
\end{figure}

\begin{figure}[H]
\centering
\scalebox{.5}{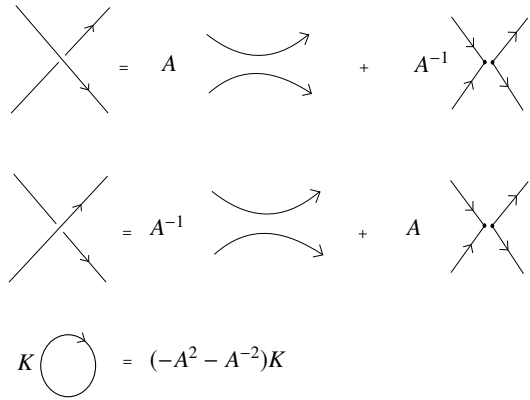}
\caption{Arrow state sum expansion}
\label{fig:arrs}
\end{figure}

\begin{figure}[H]
\centering
\scalebox{1}{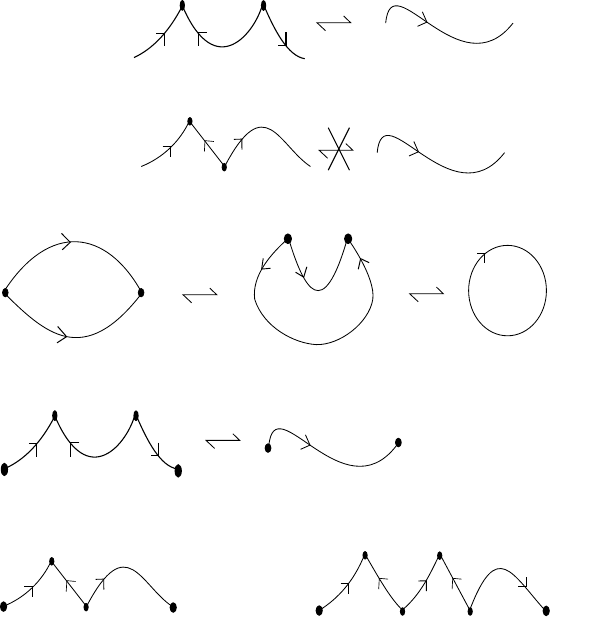}
\caption{Cusps}
\label{fig:cusp}
\end{figure}

Furthermore, the {\it affine index polynomial}, given in \cite{GK1}, is induced by a non-trivial biquandle structure on knotoid diagrams (see Figure~\ref{fig:aff}), and a number of biquandle coloring invariants were studied in \cite{GN1}.

There is also a well-defined parity assigned to crossings of (planar or spherical) knotoid diagrams. The \textit{Gaussian parity} is a mapping that assigns each crossing of a knotoid diagram to either $0$ or $1$. The importance of an existing parity for knotoids comes from the fact that there is no nontrivial parity for classical knot diagrams. Some parity based invariants for knotoids such as the \textit{odd writhe} and the \textit{parity bracket polynomial} \cite{Man} were studied in \cite{GK1}.  

Proper knotoids give rise to interesting questions related to their intrinsic nature, such as the crossing-wise distance between the endpoints, the so-called \textit{height} in \cite{GK1} (named {\it complexity} originally \cite{Tu}). More precisely, the height of a knotoid diagram in $S^2$  is the least number of crossings created when the endpoints are joined up by an underpassing strand. The  height of a knotoid is the minimum of the heights of knotoid diagrams in its equivalence class and so forms an invariant for knotoids. 

The first and the second listed authors showed that the affine index polynomial and the arrow polynomial establish lower bound estimations for the height of a knotoid. More precisely we have the following lower bound estimations for the height of a knotoid.

\begin{theorem}[\cite{GK1}] 
The height of a knotoid is greater than or equal to the maximum degree of its affine index polynomial.
\end{theorem}

\begin{theorem}[\cite{GK1}] 
The height of a knotoid is greater than or equal to the $\Lambda$-degree of its arrow polynomial.
\end{theorem}

The interested reader is referred to \cite{GK2,GKKho} for ongoing works on knotoids regarding the parity aspect of knotoids and a Khovanov homology construction in analogy with the Khovanov homology for virtual knots/links, respectively.

\section{The theory of braidoids}\label{sec:braidoids}

In this section we review the fundamental notions of braidoids introduced by the first and the last listed authors \cite{GL1,GL2}. Braidoids are defined so as to form a braided counterpart theory to the theory of planar knotoids.

\subsection{Braidoid diagrams}
Let $I$ denote the unit interval $[0,1] \subset \mathbb{R}$. A \textit{braidoid diagram} $B$ is an immersion of a finite union of arcs into $I \times I\subset \mathbb{R}^2$. The images of arcs are called \textit{strands} of $B$. There are only finitely many intersection points among the strands of $B$  which are  transversal double points endowed with over/under data, the crossings of $B$. We identify $\mathbb{R}^2$ with the $xt$-plane with the $t$-axis directed downward.  Following the  orientation of $I$, each strand is naturally oriented downward, with no local maxima or minima, so that it intersects a horizontal line at most once. 

\begin{figure}[H]
\centering \includegraphics[width=.9\textwidth]{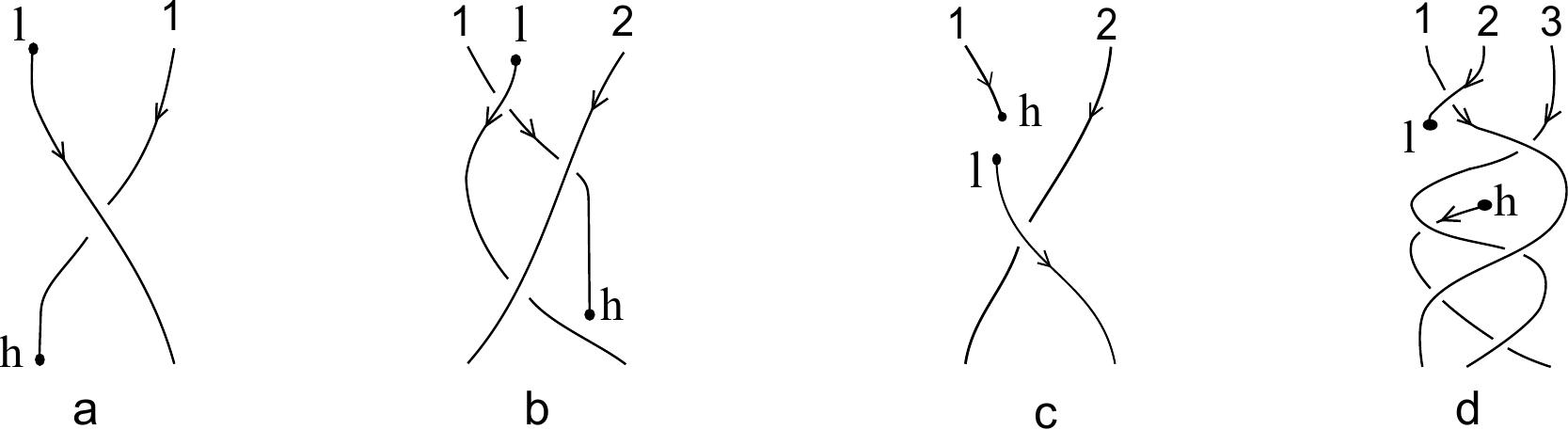} 
\caption{Some examples of braidoid diagrams}
\label{fig:basicbrds}
\end{figure}

 A braidoid diagram has two types of strands, the classical strands and the free strands. A \textit{classical strand} is as a braid strand with two ends, one lying  in $I\times\{0\}$ and the other lying  in $I\times\{1\}$.  A \textit{free strand} is a strand that either has  one of its ends lying  in $I\times\{0\}$ or in $I\times\{1\}$ and the other end lying anywhere in $I \times I$  or it has two of its ends lying anywhere in $I\times I$. There are two such free ends, called the \textit{endpoints} of $B$ and are denoted by a vertex to be distinguished from the fixed ends. For examples see Figure~\ref{fig:basicbrds}. The two endpoints are called the \textit{leg} and the \textit{head} and are denoted by $l$ and $h$ respectively, in analogy with the endpoints of a knotoid diagram. The head is the endpoint that is terminal for a free strand while the leg is the starting endpoint for a free strand, with respect to the orientation.  

The other ends of the strands of $B$ are named \textit{braidoid ends}. Each braidoid end is numbered accordingly to its order from left to right. Braidoid ends lie equidistantly and two braidoid ends having the same order on $\{ t=0 \}$ and $\{ t=1 \}$ are vertically aligned. Two braidoid ends whose orders coincide are called \textit{corresponding ends}. See the examples in Figure~\ref{fig:basicbrds}. 
\subsection{Isotopy moves on braidoid diagrams}
 We allow on braidoid diagrams the oriented Reidemeister moves $\Omega_2$ and $\Omega_3$ (recall Figure~\ref{fig:om}), which preserve the downward orientation of the arcs in the move patterns. In addition to these moves, the endpoints of a braidoid diagram can be pulled up or down in the vertical direction by a \textit{vertical move}, and right or left in the horizontal direction by a \textit{swing move} in the vertical strip determined by the neighboring corresponding ends, as long as they do not intersect or cross through any strand of the diagram. That is, the pulling of the leg or the head over or under any strand is forbidden. It is clear that allowing both  forbidden moves cancels any braiding of the free strands. 

The {\it braidoid isotopy} is induced by keeping the braidoid ends fixed on the top and bottom lines ($t=0$ and $t=1$, respectively) but allowing the Reidemeister moves $\Omega_2$ and $\Omega_3$ and planar $\Delta$-moves, as well as the swing and the vertical moves for the endpoints. Two braidoid diagrams are isotopic if one can be obtained from the other by a finite sequence of the above moves. An equivalence class of isotopic braidoid diagrams is a {\it braidoid}.

\subsection{A closure on braidoids}\label{sec:closure}

One way to define a closure operation on braidoid diagrams in order to obtain  planar (multi)-knotoid diagrams is by adding an extra property to braidoid diagrams. More precisely, each pair of the corresponding ends in a braidoid diagram is labeled either $o$ or $u$, standing for `over' or `under', respectively. We attach the labels next to the braidoid ends lying at the top line and call the diagram a \textit{labeled braidoid diagram}. Two labeled braidoid diagrams are called \textit{isotopic} if their braidoid ends possess the same labeling and they are isotopic as unlabeled diagrams. The corresponding equivalence classes are called \textit{labeled braidoids}. 

\begin{figure}[H]
\centering 
\includegraphics[width=.23\textwidth]{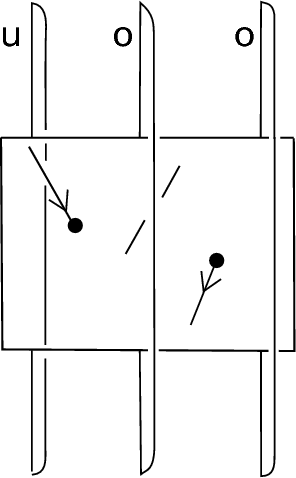}
\caption{The closure of an abstract labeled braidoid diagram}
\label{fig:c}
\end{figure}

Let $B$ be a labeled braidoid diagram. The \textit{closure} of $B$, denoted $\widehat{B}$, is a planar (multi)-knotoid diagram that results from $B$ by the following  topological operation: each pair of corresponding braidoid ends of $B$ is joined up  with a straight arc (with slightly tilted extremes) that lies on the right of the line of the corresponding braidoid ends in a distance arbitrarily close to this line so that none of the endpoints is located between the line and the closing arc. The closing arc goes entirely over or entirely under the rest of the diagram according to the label of the ends.  See Figure~\ref{fig:c} for an abstract illustration of the closure of a labeled braidoid diagram. 

\begin{proposition}[\cite{GL1,GL2}] 
The closure operation induces a well-defined map from the set of labeled braidoids to the set of planar (multi)-knotoids.
\end{proposition}



\subsection{How to turn a knotoid into a braidoid?} \label{sec:braiding}

 J.W. Alexander proved in 1923 that any oriented classical knot/link can be represented by an isotopic knot/link diagram in braided form  \cite{Al}.  The proof of the Alexander theorem by the last listed author \cite{Lathesis,LR1} utilizes the \textit{$L$-braiding moves}. In \cite{GL1,GL2} the first and the last listed authors  proved the following analogue of the Alexander theorem for (multi)-knotoids by utilizing these moves.

\begin{theorem}[\cite{GL1, GL2}]\label{thm:alex}
Any  (multi)-knotoid diagram in $\mathbb{R}^2$ is isotopic to the closure of a labeled braidoid diagram.
\end{theorem}

Let $K$ be a (multi)-knotoid diagram whose plane is equipped with the top-to-bottom direction. The basic idea for turning $K$ into a braidoid diagram is to keep the arcs that are oriented downward, with respect to the top-to-bottom direction, and to eliminate the ones that are oriented upward (\textit{up-arcs}), producing at the same time pairs of corresponding braidoid strands, such that the (multi)-knotoid resulting after closure is isotopic to $K$. The elimination of the up-arcs is done by the braidoiding moves.

\begin{figure}[H]
\centering
\scalebox{.5}{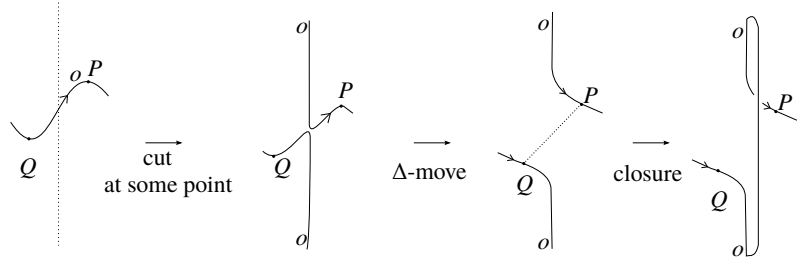}
\caption{An $L$-braidoiding move}
\label{fig:brdngone}
\end{figure}

An \text{$L$-braidoiding move} consists in cutting an up-arc at a point and pulling the resulting two pieces, the upper upward to the line $t=1$  and the lower downward to the line $t=0$, both entirely {\it over} or {\it under} the rest of the diagram. The resulting pieces are pulled so that their ends are kept aligned vertically with the cut-point. Finally the lower piece is turned into a braidoid strand by $\Delta$-moves. See Figure~\ref{fig:brdngone}. For the purpose of closure, the resulting pair of strands is labeled $o$ or $u$ depending on our choice we make for pulling the upper and lower pieces during the braidoiding move. The reader is referred to Figure~\ref{fig:brdngone} for an illustration of an $L$-braidoiding move where it can be also verified that the closure of the resulting strands labeled $o$ is isotopic to the initial up-arc.

\begin{figure}[H]
\centering
\includegraphics[width=.8\textwidth]{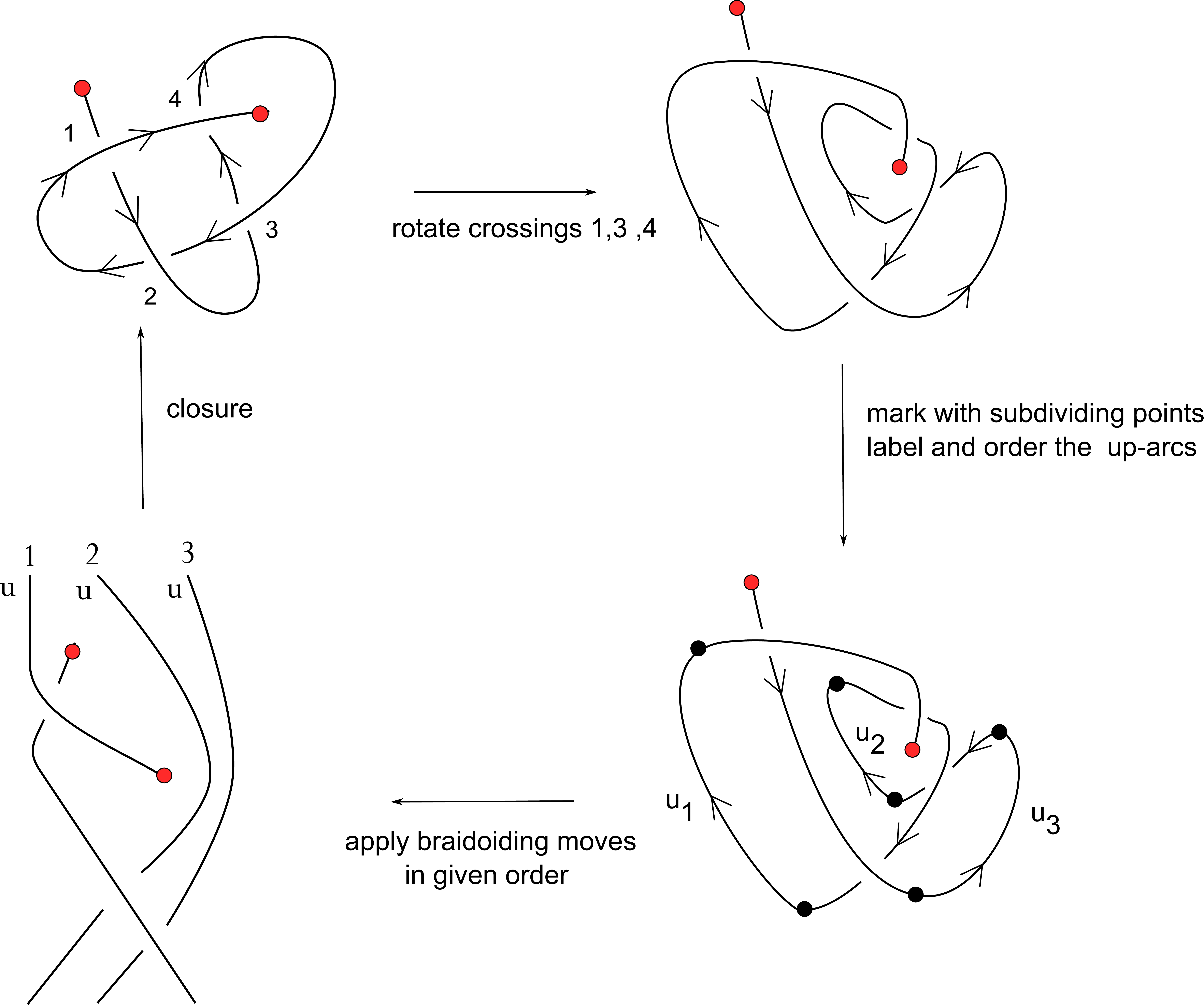}
\caption{An example showing how the algorithm works}
\label{fig:alg}
\end{figure}

Theorem~\ref{thm:alex} is proved by applying any one of the braidoiding algorithms. These algorithms are both based on the $L$-braidoiding moves. In Figure~\ref{fig:alg} we illustrate the steps of the algorithm that is given in \cite{GL2}. As for the algorithm in \cite{KL1} turning any virtual knot/link diagram into a virtual braid diagram, we start by rotating each crossing containing one or two up-arcs by $\frac{\pi}{2}$ or $\pi$, respectively, so that we end up with a knotoid diagram whose up-arcs are all free of any crossings. All of these `free' up-arcs are given an order and a labeling of {\it o} or {\it u} each, and are eliminated by the $L$-braidoiding moves one by one. The algorithm terminates in finite steps and results in a labeled braidoid diagram in Figure~\ref{fig:alg}. The algorithm in \cite{GL2} which is based on \cite{LR1}, uses braidoiding moves for up-arcs in crossings and is more appropriate for establishing Markov-type theorems for braidoids (see Theorem~\ref{thm:markov}). Yet, an added value of the first algorithm is the following consequence: \textit{Any (multi)-knotoid diagram is isotopic to the uniform closure of a braidoid diagram} \cite{Gthesis,GL2}.

\subsection{$L$-equivalence}\label{sec:l}

It is clear that due to the choices made in order to prepare a (multi)-knotoid diagram for a braidoiding algorithm (such as subdivision and labeling of the up-arcs) as well as knotoid isotopy moves, we obtain different braidoid diagrams with possibly different numbers of strands and labels. The  question that would lead to a Markov-like theorem for braidoids is to ask how these braidoid diagrams are related to each other. Clearly, the braidoid isotopy does not change the number of strands nor the labeling, so braidoid isotopy is not sufficient for determining such relations. The first and the last listed authors showed in \cite{GL1, GL2} that the $L$-moves on braidoid diagrams provide an answer to this question.

\begin{figure}[H]
\includegraphics[width=.6\textwidth]{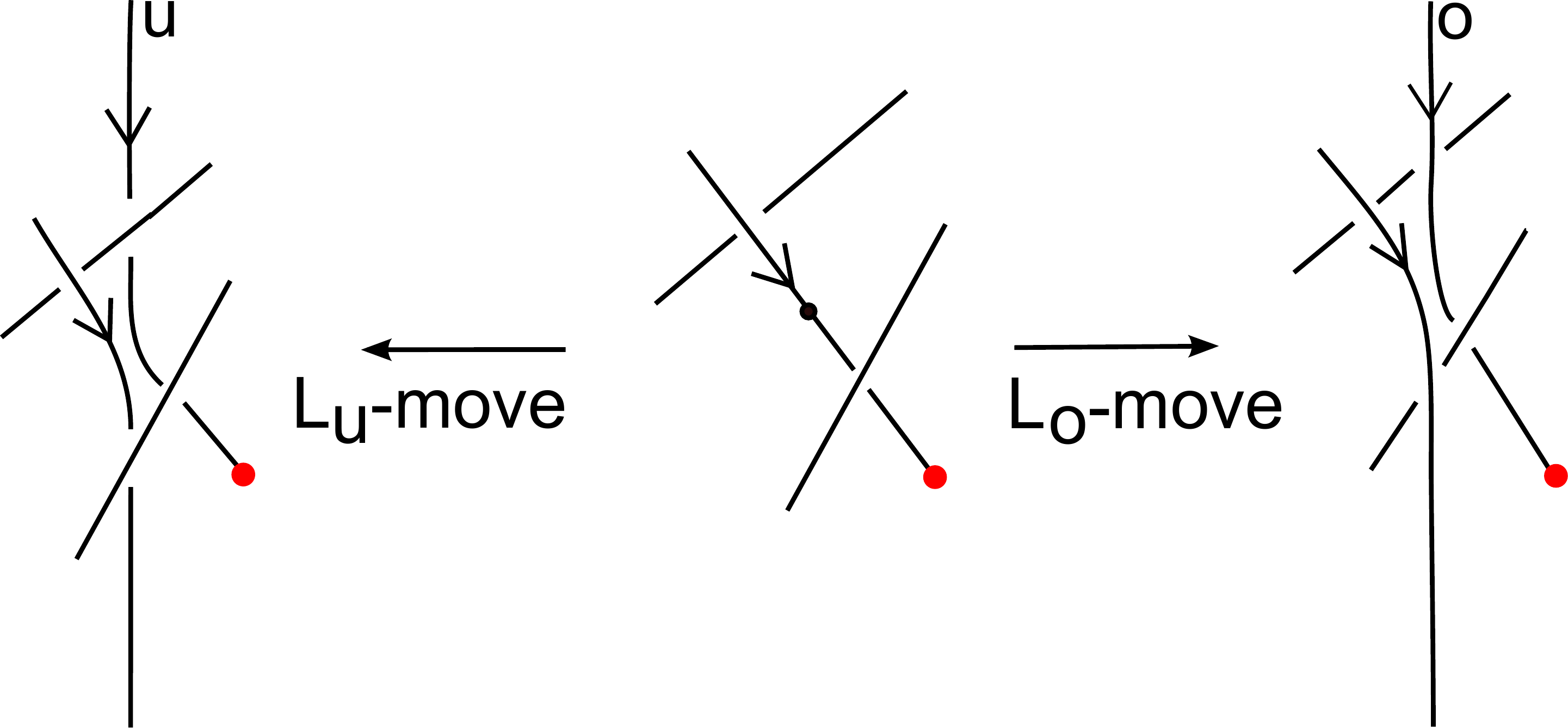}
\caption{$L$-moves}
\label{fig:alfa}
\end{figure}

An {\it $L$-move} on a braidoid diagram $B$ consists in cutting a strand of $B$ at an interior point, not  vertically aligned with a braidoid end or an endpoint or a crossing, and then pulling the resulting ends away from the cutpoint to the top and bottom of $B$ respectively, keeping them aligned with the cutpoint, and so as to create a new pair of corresponding braidoid strands. See Figure~\ref{fig:alfa} for an illustration of an $L$-move.
There are two types of $L$-moves, namely $L_o$ and $L_u$. For an $L_o$-move the pulling of the resulting new strands is entirely {\it over} the rest of the diagram. For an $L_u$-move  the pulling of the new strands is entirely {\it under} the rest of the diagram. The two  resulting strands are both labelled according to the type of the $L$-move.  See Figure~\ref{fig:alfa}.

The above definition provides us with the following result, which is an analogue of the classical Markov theorem.  

\begin{theorem}[\cite{Gthesis,GL2}]\label{thm:markov}
The closures of two labeled braidoid diagrams are isotopic (multi)-knotoids in $\mathbb{R}^2$ if and only if the labeled braidoid diagrams are related to each other by a finite sequence of $L$-moves and braidoid isotopy moves.
\end{theorem}

\section{Applications}\label{sec:app}

\subsection{Applications to the study of proteins}

\begin{figure}[H]
\includegraphics[width=.7\columnwidth]{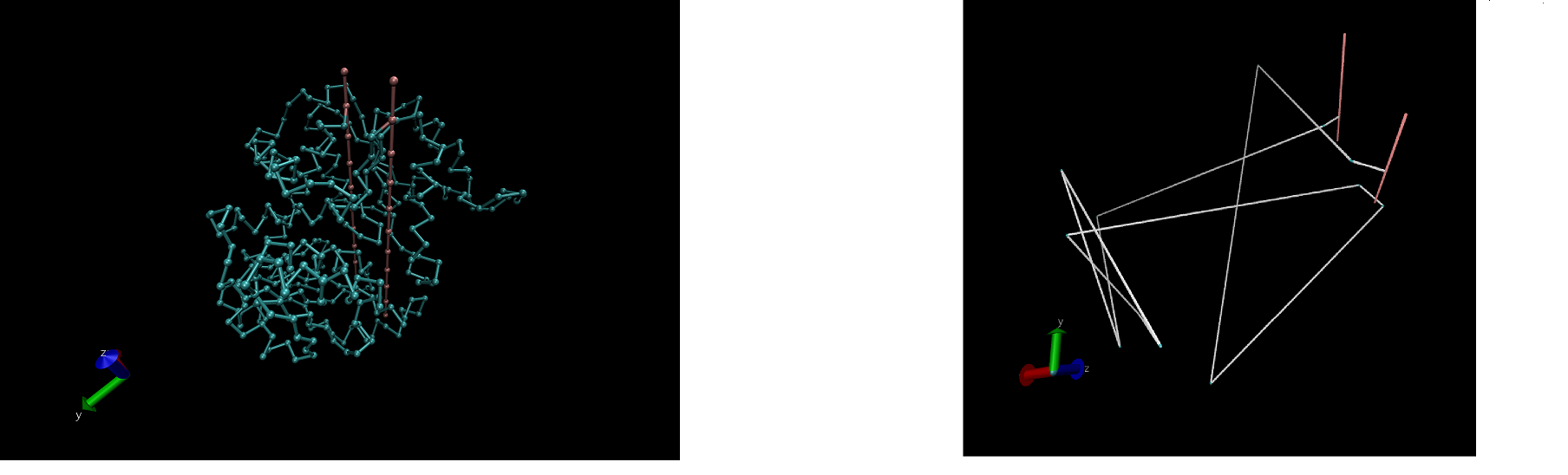}
\caption{The configuration of the backbone of the protein 3KZN in 3D and its simplified configuration; image from \small{Dimos Goundaroulis, private communication} }
\label{fig:protein}
\end{figure}

\begin{figure}[H]
\includegraphics[width=.7\columnwidth]{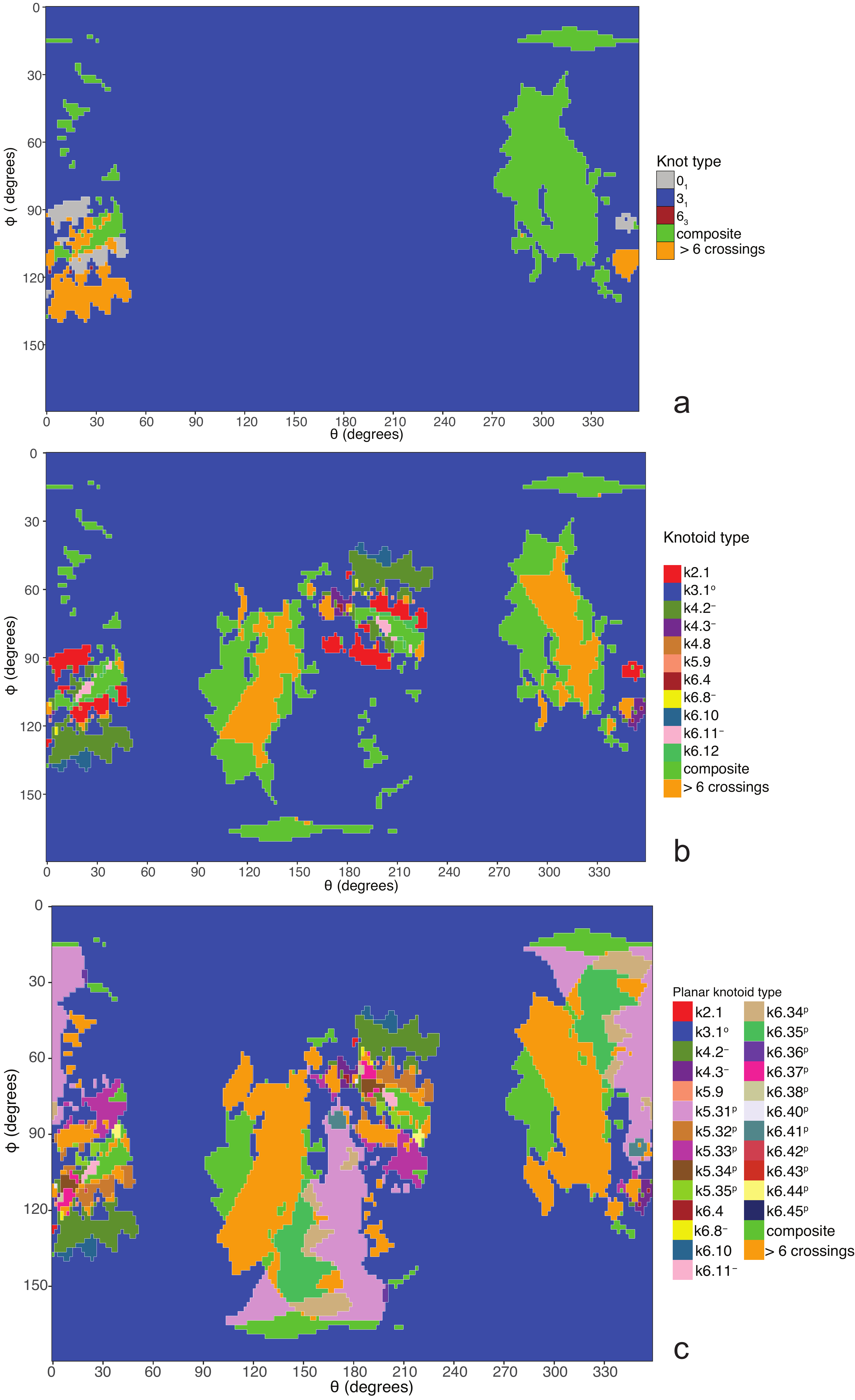}
\caption{Atlases showing the topological analysis of 3KZN via knots, spherical knotoids and planar knotoids, image from \cite{GGLDSK}}
\label{fig:atlas}
\end{figure}

The correspondence of line isotopy classes of open space curves and isotopy classes of planar knotoids suggests that a topological analysis of linear physical structures lying in $3$-dimensional space can be done by simulating them by open space curves and by taking their orthogonal projections.
 Through this idea, knotoids, both in $S^2$ and $\mathbb{R}^2$, have found important applications in the study of open protein chains \cite{GDBS} (see Figure~\ref{fig:protein} for an example), as well as of open protein chains with chemical bonds via introducing the notion of  \textit{bonded knotoids} \cite{GGLDSK}. 
In these papers, open protein chains are studied via their projections into planes. The corresponding knotoid classes are considered in the 2-sphere and in the plane, and they are classified by using the (normalized) bracket polynomial and the Turaev's loop bracket polynomial \cite{Tu, GGLDSK}, respectively. In Figure~\ref{fig:atlas} we see three atlases that contain colored regions. Each of these colored regions corresponds to one topological class of the protein 3KZN when it is closed to some knot, and when it is projected to a knotoid, a spherical knotoid and a planar knotoid, respectively. As seen from Figure~\ref{fig:atlas}, the number of colors increases as we go from the knot representation to the planar knotoid representation. By this data, the authors concluded that planar knotoids yield a more refined analysis for understanding the topological structure of open protein chains than knots or spherical knotoids \cite{GGLDSK}. This is due to the facts that more knotoids close to the same knot  and that the classification of knotoids in the plane is more refined than the classification of spherical knotoids. For example a trivial knotoid in $S^2$ may happen to be a non-trivial knotoid when considered in $\mathbb{R}^2$ as we discussed in Section~\ref{sec:knotoid}. See also \cite{DRGDSMRSS} on the subject.

\subsection{Elementary blocks and a proposed application of braidoids}\label{sec:elm}

\begin{figure}[H]
\includegraphics[width=1\columnwidth]{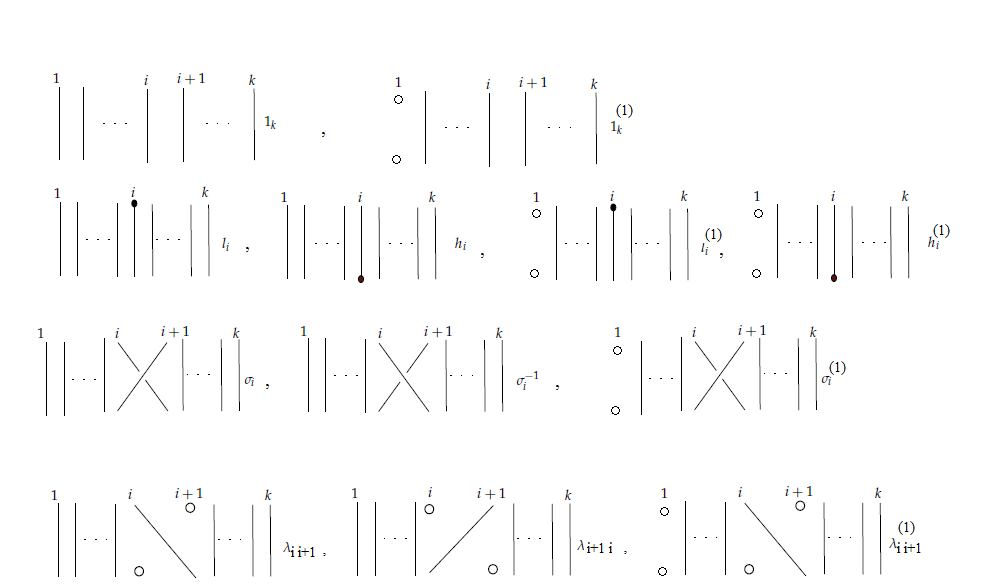}
\caption{Elementary $k$-blocks}
\label{fig:elt}
\end{figure}

As demonstrated in \cite{GL1,Gu1}, any braidoid diagram can be divided into a finite number of horizontal stripes, each containing one of the blocks that are depicted in Figure~\ref{fig:elt}. The blocks consist of the classical braid generators, the identity elements containing one endpoint, along with their extensions by the \textit{implicit points}, which are empty nodes put along the vertical direction of the endpoints, and along with the \textit{shifting blocks}, which result from the change of positioning of strands before or after the  appearance of an endpoint.  
A product on the set of elementary $n$-blocks and relations with respect to this  product, induced by the isotopy moves of the braidoid diagrams, is further explored in \cite{Gu1}. Then, any braidoid diagram on $n$ strands can be read, from top to bottom, as a word that corresponds to a combination of finitely many elementary $n$- or $(n+1)$-blocks. 

\begin{figure}[H]
\centering \scalebox{.3}{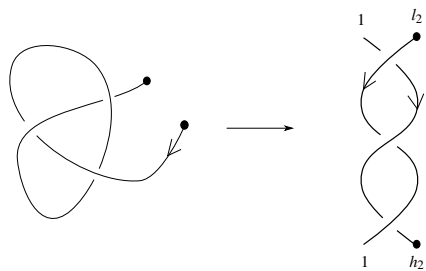}
\caption{The knotoid of the protein 3KZN and a corresponding  braidoid diagram with algebraic expression $l_2\sigma_1^3h_2$}
\label{fig:prot}
\end{figure}

Any planar knotoid diagram can be turned into a (labeled) braidoid diagram by  Theorem~\ref{thm:alex}, so it can be represented by an expression in terms of elementary blocks. This suggests an algebraic encoding for open protein chains or, in general, for linear polymer chains: they can be projected to  planes and the resulting knotoid diagrams can be turned into braidoid diagrams that have algebraic expressions. An example is illustrated in Figure~\ref{fig:prot}, where the knotoid corresponding to  protein 3KZN is turned into a braidoid diagram, which is represented by the word $l_2\sigma_1^3h_2$ in elementary blocks.

\begin{acknowledgement}
The research of Sofia Lambropoulou and Neslihan G\"ug\"umc\"u has been co-financed by the European Union (European Social Fund - ESF) and Greek national funds through the Operational Program ``Education and Lifelong Learning'' of the National Strategic Reference Framework (NSRF) - Research Funding Program: THALES: Reinforcement of the interdisciplinary and/or inter-institutional research and innovation, MIS: 380154. Louis H. Kauffman' s work was supported by the Laboratory of Topology and Dynamics, Novosibirsk State University (contract no. 14.Y26.31.0025 with the Ministry of Education and Science of the Russian Federation)
 \end{acknowledgement}






\end{document}